\documentclass[10pt]{article}

\usepackage{graphicx}
\usepackage{amsmath,amssymb,amsfonts,color}

\newtheorem{lemma}{Lemma}

\newtheorem{proposition}{Proposition}

\renewcommand{\P}{\mathbb{P}}
\newcommand{\E}{\mathbb{ E}}
\newcommand{\I}{{\bf I}}

\newcommand{\fd}{\hspace{2mm} $\square$ \vspace{1mm}}
\font\simb=msam10
\newcommand{\trn}{\vspace{1.5mm}{\simb I}\ \rm }

\begin{document}

\title{Almost sure convergence of Polya urn schemes}

\author{Ricardo V\'elez\\
{\small\it Statistics Deparment, UNED, Spain}}
\date{}

\maketitle

\begin{abstract}
For the most general Polya urn schemes, we establish the almost sure convergence of its composition.
The only requirement is that there are always enough balls of both colors, so
that the extractions can be indefinitely  pursued according to the specifications of the model.
We also consider the method for determining the probability of fulfilling this requirement, as a
function of the initial number of balls of each color.
\end{abstract}

\section{Polya's Urn schemes}

An urn scheme is concerned with an urn containing balls of two different colors:
$\mathcal{A}$mber and $\mathcal{B}$lue.
Successive random extractions are performed  and the  ball that has been drawn is returned
to the urn  together with a certain number of  additional balls of each color.
Initially there is a total of $t_0$ balls in the urn,
$\alpha_0$ of which are amber and $\beta_0=t_0-\alpha_0$ are blue.
The urn scheme is characterized further by the parameters:
$$[\hspace*{2mm}\underbrace{a'\hspace*{1mm},\hspace*{2mm}a}_{amber}\hspace*{2mm};
\hspace*{2mm}\underbrace{b\hspace*{2mm},\hspace*{2mm}b'}_{blue}\hspace*{2mm}]$$
where $a'$ and $a$ stand respectively for the number of amber and blue balls introduced in
the urn when an amber ball is obtained, while $b$ amber and $b'$ blue balls are adjoined
if a blue ball is extracted.
Negative values for some of these values are usually admitted and adding  $-x$ balls means
removing $x$ balls from the urn. However, in this case it may happen that the balls of some
color are exhausted and the extractions cannot be carried on according to the rules of the model.

The history of Polya's urn schemes is reported in references \cite{JK}, \cite{KB}. A more
recent survey on this subject is the monograph \cite{Mahmoud}. The reference \cite{FDP} also
contains a review of the main papers in the area.

In \cite{VP} we have analyzed the long term behavior of Polya's urm schemes under the simplified
hypothesis that $a'+a=b+b'=c>0$. So the same number of balls is added after each extraction,
whatever the color of the obtained ball can be,  and therefore the total number $t_n=t_0+cn$
of balls in the urn grows deterministically.

Our main conclusion for this case is that the proportion $p_n$
of amber balls in the urn converges almost surely to $b/(a+b)$, whenever both colors
survive indefinitely in the urn (and, in particular, if $0\leq a,b\leq c$). For $a<0$
the blue balls are exhausted and the amber balls deplete if $b<0$. The more exceptional
case when $a=b=0$ is the well known Polya--Eggenberger model, for which
$p_\infty=\lim p_n$ exists almost surely, but is a random variable with distribution
beta$(\alpha_0/c,\beta _0/c)$.

These results hold true even if $a$ and $b$ have random values. More precisely,
assume that $\{a_k\}_{k\geq 1}$ and $\{b_k\}_{k\geq1}$ are two sequences of independent
identically distributed random variables with  finite ranges of integer values, independent
also of the sequence of colors observed. At the $k$-th stage, $a_k$ blue and $c-a_k$ amber balls are
added  when an amber ball is obtained, while $b_k$ amber and $c-b_k$ blue are introduced if a blue
ball appears. For this ``random'' model, $\lim p_n=\bar{b}/(\bar{a}+\bar{b})$ where $\bar{a}$
and $\bar{b}$ are the means of $a_k$ and $b_k$ respectively. An urn scheme will be called
non random if $a$ and $b$ have fixed values.

Now, our aim here is to consider the case when $A=a'+a$ and $B=b+b'$ are different,  so
that the total number of balls in the urn evolves randomly, since it depends on the successive
colors of the extracted balls. Specifically, we assume that $A,B>0$ in order to have an increasing
total number of balls in the urn.  Since both colors may be interchanged, we can assume without loss
of generality that $\Delta=A-B>0$. Summarizing, the more general model to be considered now is
$$[\hspace*{2mm}\underbrace{A-a_k\hspace*{1mm},\hspace*{2mm}a_k}_{amber}\hspace*{2mm};
\hspace*{2mm}\underbrace{b_k\hspace*{2mm},\hspace*{2mm}B-b_k}_{blue}\hspace*{2mm}]$$
at stage  $k$, where $A>B>0$ are given integers. We will use the following notation:
\begin{itemize}
\item[\trn] $Y_k$, the color of the ball obtained in the trial $k\geq 1$: \\
\rule{0mm}{5mm} $Y_k=1$ if it is amber,\hspace*{3mm} $Y_k=0$ if it is blue.
\item[\trn] $X_n=\sum_{k=1}^{n}Y_k$, the total number of amber balls obtained in the
first $n$ extractions.
\item[\trn] $t_n=t_0+AX_n+B(n-X_n)=t_0+nB+\Delta X_n$\\  \rule{0mm}{5mm} is the total number of
balls after the first $n$ stages. Recurrently
\begin{equation}\label{}
t_{n+1}=t_n+AY_{n+1}+B(1-Y_{n+1})
\end{equation}
\item[\trn] $\alpha_n=\alpha_0+\sum_{k=1}^{n} (A-a_k)\, Y_k+\sum_{k=1}^{n}b_k\,(1-Y_k)$ \\ \rule{0mm}{5mm}
is the total number of amber balls when $n$ extractions have been performed. It satisfies
\begin{equation}\label{}
\alpha_{n+1}=\alpha_n+(A-a_{n+1})Y_{n+1}+b_{n+1}(1-Y_{n+1}).
\end{equation}
\item[\trn] $\beta_n=t_n-\alpha_n$ \quad is the total number of blue balls after the first $n$ steps.

\item[\trn] $p_n=\alpha_n/t_n$ \quad is the proportion of amber balls in the urn after the first
$n$ stages.
\end{itemize}

Recall that $\{a_k\}_{k\geq 1}$ and $\{b_k\}_{k\geq1}$ are sequences of independent identically
distributed random variables with finite ranges of integer values, independent also of
$\{Y_k\}_{k\geq1}$. The mean values of their respective distributions will be represented by
$\bar{a}$ and $\bar{b}$.

When the ranges of $a_k$ or $b_k$ include values that are negative or such that $A-a_k$ or $B-b_k$
are negative, it is mandatory to consider
$$\tau=\inf\{n\geq1 \,| \; \alpha_n<0 \text{ or } \alpha_n>t_n\}=\inf\{n\geq1\, |\; p_n\notin [0,1]\}$$
which represents the first time when the balls of some color exhaust, so that the extractions cannot be
pursued according to the model specifications. Of course $\{\tau=\infty\}$  is a sure event if
$0\leq a_k\leq A$ and $0\leq b_k\leq B$; but otherwise it can be $\P\{\tau<\infty\}\in (0,1)$.
In fact, assuming that $A-a_k,a_k,b_k,B-b_k$  have only small negative values,
$\P\{\tau=\infty\}$ quickly increases as $\alpha_0,\beta_0\nearrow\infty$.

On $\{\tau >n\}$ the $n$-th extraction can be performed and $p_n$ is the probability of getting an
amber ball.
For trajectories in $\{\tau=\infty\}$, all the sequence $\{p_n\}$ is defined and we want to
analyze its limit behavior by means of martingale arguments. This is done in section 3 after some preliminary results.

\section{Preliminary results}

For the random Polya's urn scheme $[A-a_k ,a_k\hspace*{1mm};\hspace*{1mm}b_k,B-b_k]$,
let ${\cal F}_n=\sigma(\{Y_k,a_k,b_k\}_{1\leq k\leq n})$ be the $\sigma$-field of events
depending on the results of the $n$ first extractions.

It is a natural guess that the asymptotic number of amber balls obtained equals the
limit proportion of amber balls in the urn. Here is the precise result.

\begin{lemma}
In $\{\tau =\infty\}$, if $p_n\rightarrow p_\infty$ almost surely then also
$X_n/n\longrightarrow p_\infty$ almost surely.
\end{lemma}
{\em Proof:}  Let us consider  $Z_k=Y_k-p_{k-1}$, for which $\E[Z_k|{\cal F}_{k-1}]=0$.
Then we have
$$\E[Z_k^2]=\E[Y_k]-2\E[Y_kp_{k-1}]+\E[p_{k-1}^2]=\E[p_{k-1}]-\E[p_{k-1}^2]\leq \frac{1}{4}$$
and therefore $\sum_{k=1}^{\infty}k^{-2}\E[Z_k^2]<\infty$. Thus, the strong law of large numbers for
martingales given in \cite[VII, Theorem 3]{Feller} asserts that, almost surely in $\{\tau =\infty\}$,
$$\frac{1}{n}\sum_{k=1}^{n}Z_k=\frac{X_n}{n}-\frac{1}{n}\sum_{k=1}^{n}p_k\longrightarrow 0
\hspace*{7mm} \text{or} \hspace*{7mm}   \frac{X_n}{n}\longrightarrow p_\infty.\hspace*{6mm} \square$$

Although we will prove that ${p_n}$ behaves as a sub or supermartingale,
being general results, the martingale convergence theorems do not give any information
about the limit. In the present setting the following result holds.

\begin{lemma}
If $p_n\rightarrow p_\infty$ almost surely in $\{\tau=\infty\}$, then $p_\infty$ takes
values between the roots of the polynomial
\begin{equation}\label{}
\omega(x)=\Delta x^2-(\Delta-\bar{a}-\bar{b})x-\bar{b}.
\end{equation}
\end{lemma}
{\em Proof:} First assume that $X_n\nearrow\infty$ and therefore
$$\frac{1}{X_n}\sum_{k=1}^{n}Y_ka_k\rightarrow\bar{a}\hspace*{5mm}\text{and}\hspace*{5mm}
\frac{1}{X_n}\sum_{k=1}^{n}Y_kb_k\rightarrow \bar{b}$$
whatever values the sequence $\{Y_k\}_{k\geq1}$ may have. Now observe that
$$p_n\frac{t_n}{n}=\frac{\alpha_n}{n}=\frac{\alpha_0}{n}+A\frac{X_n}{n}-
\frac{X_n}{n}\frac{1}{X_n}\sum_{k=1}^{n}a_kY_k+\frac{1}{n}\sum_{k=1}^{n}b_k
-\frac{X_n}{n}\frac{1}{X_n}\sum_{k=1}^{n}b_kY_k.$$
and $t_n/n\rightarrow B+\Delta p_\infty$.
Thus any limit value $p_\infty$ must satisfy the equation
\begin{equation}\label{ecp_inf}
Bp_\infty+\Delta p_\infty^2=(A-\bar{a}-\bar{b})p_\infty+\bar{b}.
\end{equation}
Therefore $p_\infty$ can take the value zero (if $X_n$ is bounded) or some
root of the polynomial $\omega$.

But, the Borel--Cantelli lemma in \cite[Theorem 5.3.2]{D} gives
$\{X_n\rightarrow\infty\}=\{\sum_{n=1}^{\infty}p_n=\infty\}$,
so that, if $X_n$ remains bounded, it is $\sum_{n=1}^{\infty }p_n<\infty$;
and, since $p_n=0$ or $p_n\geq1/(t_0+nB)$, it must be $p_n=0$ for all large enough $n$.
Moreover $p_n=0=p_{n+1}$ implies that $b$ has a non random vanishing value and
$0=\lim p_n$ is a root of $\omega(x)$. \fd

\section{Martingale analysis}

Assuming that $\tau >n+1$, it is
\begin{align*}
p_{n+1}&=Y_{n+1}\frac{\alpha_n+A-a_{n+1}}{t_n+A}+(1-Y_{n+1})\frac{\alpha_n+b_{n+1}}{t_n+B}
\end{align*}
so that
\begin{align*}\label{}
\E[p_{n+1}|{\cal F}_n]&=p_n\frac{\alpha_n+A-\bar{a}}{t_n+A}+(1-p_n)\frac{\alpha _n+\bar{b}}{t_n+B}\\
&=p_n\frac{p_n+(A-\bar{a})/t_n}{1+A/t_n}+(1-p_n)\frac{p_n+\bar{b}/t_n}{1+B/t_n}
\end{align*}
and
\begin{equation}\label{}
\E[p_{n+1}|{\cal F}_n]-p_n=\frac{h(p_n,t_n)}{(1+A/t_n)(1+B/t_n)}
\end{equation}
where (recall that $\Delta=A-B$)
\begin{equation}\label{}
h(p,t)=-p^2\;\frac{\Delta}{t}+p\;\frac{\Delta-\bar{a}(1+B/t)-\bar{b}(1+A/t)}{t}+\frac{\bar{b}(1+A/t)}{t}
\end{equation}
is a convex parabolic function of $p$ (\footnote{For $\Delta=0$, $h(p,t)$ is a linear function vanishing at
$\bar{b}/(\bar{a}+\bar{b})$. This remark suffices to obtain the conclusions of \cite{VP} by means of simple
arguments similar to those to be presented later.}).
Observe that $$h(0,t)=\bar{b}/t(1+A/t)\hspace*{6mm}\text{and} \hspace*{6mm}h(1,t)=-\bar{a}/t(1+B/t).$$
One can now distinguish various cases.

\subsection{ The case $\bar{a}\leq0\leq\bar{b}$}

For $\bar{a}\leq0\leq\bar{b}$, it is $h(p,t)\geq0\;\forall p\in[0,1]$. Thus, in $\{\tau=\infty \}$,
$\{p_n\}_{n\geq 1}$ is a bounded submartingale with respect to $\{{\cal F}_n\}$ and
converges to a random variable $p_\infty\in[0,1]$ almost surely and in $L^1$.
Moreover, Lemma 2 allows to conclude
\begin{proposition}
\begin{itemize}
\item[i)] If $\bar{a}<0\leq\bar{b}$  then  $\P\{\tau<\infty\}=1$.
\item[ii)] For $\bar{a}=0\leq\bar{b}$ it is $\P(\{\tau <\infty\}\cup\{p_\infty=1\})=1$.
\end{itemize}
\end{proposition}

In fact, for $\bar{a}<0<\bar{b}$, $\omega (x)$ takes only negative values in $[0,1]$;
thus $p_\infty$ cannot exist and the conclusion follows. For $\bar{a}<0=\bar{b}$, $\omega(x)$
has only the root $x=0$, but the possibility of being $p_\infty=0$ is excluded since
this  would imply $\E[p_n\,\I_{\{\tau=\infty\}}]\nearrow 0$. Thus $\P\{\tau=\infty\}=0$ follows.

If $\bar{a}=0<\bar{b}$ the only non negative root of $\omega(x)$ is $x=1$ and the result is proved.

For $\bar{a}=0=\bar{b}$, $p_\infty$ can take also the value 0; but
$\E[p_n\I_{\{p_\infty=0\}}]\nearrow 0$ implies $\P\{p_\infty=0\}=0$.

\vspace{2mm}
This last case, with $a$ and $b$ non random, is an unbalanced Polya-Eggenberger model, in which
the number of accompanying amber balls greater than the number of accompanying blue balls.
The interchangeability of the variables $\{Y_k\}$
(establishing the beta distribution of $p_\infty$)
fails and only $p_\infty=1$ may happen. So the beauty of the Polya--Eggenberger model is due
to its symmetry.

\subsection{ The case $\bar{a},\bar{b}>0$}

Since $h(0,t)>0$ and $h(1,t)<0$, there is a unique root $p^{\star}_n\in(0,1)$ such that
$h(p,t_n)>0$ for $p<p^{\star}_n$ and $h(p,t_n)<0$ for $p>p^{\star}_n$. But, as $n\rightarrow\infty$,
$t_n\nearrow\infty $ and  $p_n^{\star}$ converges to
\begin{equation}\label{p*}
p^{\star}=\frac{\Delta-\bar{a}-\bar{b}+\sqrt{(\Delta-\bar{a}-\bar{b})^2+4\Delta b}}{2\Delta}
\end{equation}
which is also the only root in $(0,1)$ of $\omega(x)$. We will prove that $\{p_n\}$ converges
almost surely and Lemma 2 will give the value of the limit.
\begin{proposition}
If $\bar{a},\bar{b}>0$, in $\{\tau=\infty\}$ it is $p_n\rightarrow p^{\star}$ almost surely.
In other words $\P(\{\tau<\infty\}\cup \{p_n\rightarrow p^{\star}\})=1$.
\end{proposition}
{\em Proof:} $\{p_n\}$ is a submartingale as long as $p_n\leq p^{\star}_n$ and is a supermartingale
when $p_n\geq p^{\star}_n$.

Within the event $\{\tau =\infty\}$,  $\bar{p}=\limsup p_n$ and
$\underline{p}=\liminf p_n$ both exist and we can consider the event
$C=\{\tau=\infty,\underline{p}<\bar{p}\}$.
Since  $p^{\star}_n$ gets arbitrarily close to $p^{\star}$, $C\cap \{\bar{p}<p^{\star}\}$ has probability
zero because those trajectories are, from some $n$ onwards, trajectories of a submartingale with two
different cluster points. Similarly $C\cap\{\underline{p}>p^{\star}\}$ has probability zero.
Hence,
$C\cap\{\underline{p}<p^{\star}<\bar{p}\}$ differs from $C$ by a set of probability zero.

Now, let $C_\delta=C\cap\{\underline{p}<p^{\star}<p^{\star}+\delta<\bar{p}\}$. The trajectories in $C_\delta$
must perform an infinite number of upcrossings of the interval $(p^{\star},p^{\star}+\delta)$ through positive
steps of size
$$p_{n+1}-p_n=\left\{\begin{array}{ll} \displaystyle
\frac{A-a_{n+1}-p_n A}{t_n+A} & \text{if } Y_{n+1}=1\\ \displaystyle
\frac{b_{n+1}-p_nB}{t_n+B}   & \text{if }  Y_{n+1}=0 \rule{0mm}{6mm}
\end{array}\right.$$
that is less than any $\varepsilon>0$ for $n$ large enough. Therefore the probability of $C_\delta$ is bounded
by $[(p^{\star}+\delta) \vee (1-p^{\star})]^{\delta/\varepsilon}$ and, this being true for any $\varepsilon>0$, it
should be $\P(C_\delta)=0$. A similar reasoning shows that
$C'_\delta=C\cap\{\underline{p}<p^{\star}-\delta<p^{\star}<\bar{p}\}$ has also probability zero and consequently
$\{\tau<\infty\}\cup\{\underline{p}=p^{\star}=\bar{p}\}$ has probability one. \fd

\subsection{ The case $\bar{a},\bar{b}<0$}

Here, $h(p_n,t_n)<0$ for $p_n<p^{\star}_n$ and $h(p_n,t_n)>0$ for $p_n>p^{\star}_n$ where $p^{\star}_n$
is the unique root in $(0,1)$ of $h(t_n,p)$. However, now $p^{\star}_n$ converges to
\begin{equation}\label{p_*}
p_{\star}=\frac{\Delta-\bar{a}-\bar{b}-\sqrt{(\Delta-\bar{a}-\bar{b})^2+4\Delta b}}{2\Delta}
\end{equation}
which is again the only root in $(0,1)$ of $\omega(x)$.
The rest of the reasoning of the last section holds without any change, so that the same conclusion holds:
\begin{proposition}
If $\bar{a},\bar{b}<0$ it is $\P(\{\tau<\infty\}\cup\{p_n\rightarrow p_{\star}\})=1$.
\end{proposition}
Because $p_n$ is now ``decreasing'' below and ``increasing'' above $p^{\star}$, surely $\P\{\tau<\infty\}$
is much larger in this case that when $\bar{a},\bar{b}>0$. It may be near 1 even for large values of $\alpha_0$
and $\beta_0$.

\subsection{ The case $\bar{b}\leq 0\leq\bar{a}$}

If $\bar{b}<0<\bar{a}$, $h(p,t)$ has negative values at $p=0$ and $p=1$ and the vertex
of the parabola is located at
$$\hat{p}=\frac{1}{2}-\frac{\bar{a}(1+B/t)+\bar{b}(1+A/t)}{2\Delta}$$
which is in $[0,1]$ when $|\bar{a}(1-B/t)+\bar{b}(1+A/t)|\leq\Delta $. At $\hat{p}$ it is
$$t\,h(\hat{p},t)=b(1+A/t)+\frac{[\Delta-a(1+B/t)-b(1+A/t)]^2}{4\Delta}.$$

Thus it will be $h(p_n,t_n)<0$ for all $p_n\in[0,1]$ and $n$ large enough, under one of the conditions:
$$\text{(i)}\hspace*{2mm}|\bar{a}+\bar{b}|>\Delta\hspace*{7mm}\text{or}\hspace*{7mm} \text{(ii)}\hspace*{2mm}
|\bar{a}+\bar{b}|\leq\Delta \hspace*{2mm}\text{ \ and \ }\hspace*{2mm} [\Delta-\bar{a}-\bar{b}]^2+4\bar{b}\Delta<0. $$
In other words, under these conditions, in $\{\tau=\infty\}$, $\{p_n\}$ will become a
supermartingale when $n$ increases and therefore it will converge almost surely to a limit $p_\infty$.
But, (i) and (ii) give that $\omega(x)$ has no roots in the interval $[0,1]$, thus $p_\infty$
cannot exist and $\P\{\tau<\infty\}=1$.

Under (i) or (ii) with $\bar{a}=0$, $\{p_n\}$ is still a supermartingale for large $n$ and $\omega(1)=0$;
but as $\E[p_n\I_{\{\tau=\infty\}}]$ decreases, it cannot be $p_\infty=1$, and $\P\{\tau <\infty \}=1$ also holds.

When $\bar{b}=0$, since  $\omega(0)=0$, the supermartingale $\{p_n\}$  may converge to $p_\infty=0$.
In summary:

\begin{proposition}
Assume that $|\bar{a}+\bar{b}|>\Delta$ or $|\bar{a}+\bar{b}|\leq\Delta$ and $[\Delta-\bar{a}-\bar{b}]^2+4b\Delta< 0$,
then
\begin{itemize}
\item[(i)] if $\bar{b}<0\leq\bar{a}$, it is $\P\{\tau<\infty\}=1$.
\item[(ii)] if $\bar{b}=0<\bar{a}$, it is $\P(\{\tau<\infty\}\cup\{p_\infty=0\})=1$.
\end{itemize}
\end{proposition}

We now consider the situation when $\bar{b}<0<\bar{a}$, $|\bar{a}+\bar{b}|\leq\Delta$ and
$[\Delta -\bar{a}-\bar{b}]^2+4\bar{b}\Delta\geq0$,
so that $h(p,t_n)=0$ has two roots, $p^{(1)}_n,p^{(2)}_n\in(0,1)$  converging respectively to
$p_{\star}$ and $p^{\star}$ (given in (\ref{p_*}) and (\ref{p*}))
which are both roots of $\omega(x)$.

Since $\{p_n\}$ is a submartingale when $p_n\in[p^{(1)}_n,p^{(2)}_n]$ and a supermartingale
if $p_n\in[0,p^{(1)}_n]\cup[p^{(2)}_n,1]$, the same reasoning of section 2.2 shows that, in $\{\tau =\infty\}$,
$\underline{p}=\liminf p_n$ and $\bar{p}=\limsup p_n$ cannot belong to the same interval of the
partition $[0,p^{(1)}], [p^{(1)},p^{(2)}],[p^{(2)},1]$. Also the same argument of section 2.2
proves that it must be $\underline{p}=p_{\star}=\bar{p}$ or $\underline{p}=p^{\star}=\bar{p}$.
Thus we conclude
\begin{proposition}
Assume that $\bar{b}\leq0\leq\bar{a}$, $|\bar{a}+\bar{b}|\leq\Delta $ and $[\Delta-\bar{a}-\bar{b}]^2+4b\Delta\geq0$,
then in $\{\tau =\infty\}$ the sequence $p_n$ converges almost surely to a random variable $p_\infty$
taking one of the values $p_{\star}\ or\ p^{\star}$.
\end{proposition}

Of course, if $\bar{b}=0$ it is $p_{\star}=0$ and, for $\bar{a}=0$, it is  $p^{\star}=1$.
The distribution of $p_\infty$
is not easy to find, but the simulations show that $\P\{p_\infty=p^{\star}\}$ is much
greater than $\P\{p_\infty=p_{\star}\}$. This is a plain consequence of the fact
that $p_{\star}$ is unstable ($p_n$ ``decreases'' under $p_{\star}$ and ``increases''
above $p_{\star}$); exactly the opposite of what happens at $p^{\star}$.

\begin{figure}[h!]
	\centering
	\includegraphics[width=110mm]{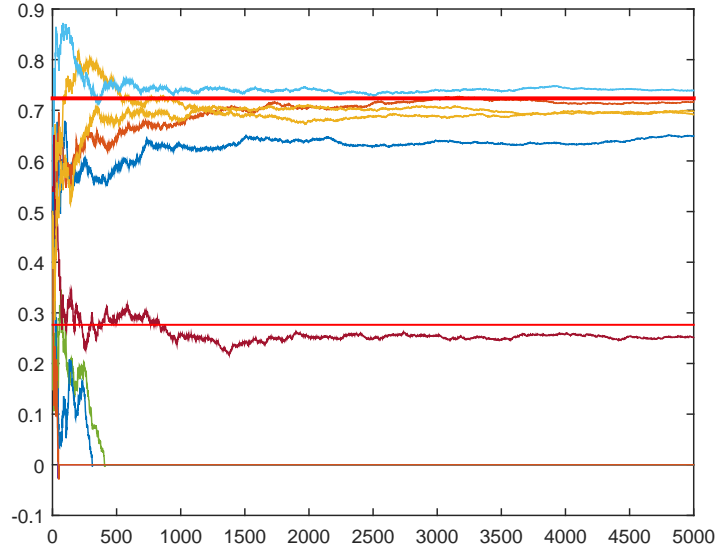}\\
	\caption{10 trajectories of the sequence $\{p_n\}$ with $A=7,B=2$ and $\bar{a}=1, \bar{b}=-1$}
	\label{fig1}
\end{figure}

A simple Matlab program (see the appendix) allows to simulate the trajectory of the sequence
$\{p_n\}$, with different values of the parameters, confirming the results of Propositions 1 to 5.
For instance, the following figure  shows 10 paths, each one corresponding to 5000 extractions of
an urn containing initially $\alpha=\beta=30$ balls of each color.
An amber ball is returned to the urn together with $A=7$ extra balls: $A-a_k$ amber and $a_k$ blue,
where each $a_k$ is chosen independently at each stage between $[-5,-2,4,7]$ with probabilities proportional to $[1,2,2,1]$
(so that  $\bar{a}=1$). Similarly the blue balls are returned to the urn with $B=2$ extra balls: $b_k$ amber and
$B-b_k$ blue, the $b_k$ being chosen in $[-5,0,4]$ with probabilities proportional to $[2,3,1]$ (and $\bar{b}=-1$).

Figure 1 is somewhat faked. The simulation has been run many times until obtaining a trajectory converging
to $p_{\star}=0'2764$. Another five converge to $p^{\star}=0'7236$ more or less slowly, while the four
remaining belong to $\{\tau<\infty\}$ with $\tau=45,50, 311,407$ all them by lack of amber balls
(trajectories in which blue balls run out are less usual with these parameters values). A very rough
estimation of the probability of the trajectories converging to $p_{\star}$ is $0'006$.

\section{About the distribution of $\tau$}

Except under the conditions $0\leq a_k\leq A$ and $0\leq b_k\leq B$ (granting that $\tau=\infty$),
it may happen that $\P\{\tau<\infty\}$ is 1 or close to 1 and therefore the proportion of converging
trajectories is very small. Then one must be interested in knowing $\P\{\tau=\infty\}$ or its
approximate value $\P\{\tau>M\}$ for large values of $M$.

To this end let
$$ q_n(t,\alpha)=\P\{\tau>M\,|\,t_n=t,\alpha_n=\alpha\}$$
for $0\leq\alpha\leq t$ and $q_n(t,\alpha)=0$ otherwise.

Assume that all $a_k$ take a value $a\in R$ with probability $r(a)$, while the
$b_k$ take a value $b\in S$ with probability $s(b)$. Then the following recurrent equation holds
\begin{equation}\label{rec}
q_n(t,\alpha)=\frac{\alpha}{t}\sum_{a\in R}r(a)q_{n+1}(t+A,\alpha+A-a)+\left(1-\frac{\alpha}{t}\right)
\sum_{b\in S} s(b)q_{n+1}(t+B,\alpha+b).
\end{equation}
As $q_M(t,\alpha)=1$ for $0\leq\alpha\leq t$, equation (\ref{rec}) may be solved backwards
so as to get $q_0(t,\alpha)$ for each $\alpha\in[0,t]$ and $t\in[t_1,t_2]$.

Such a calculation is made by the last Matlab program in the appendix, where one can fix for instance $M=t_2+fA$
with $f=100,200$, etc. However the program runs very slowly since it must get successively half a matrix of
dimension $(t_2+nA)\times (t_2+nA+1)$ for $n=M-1,M-2,\ldots,1,0$.

As a sample of the results, with the same parameters used in the preceding figure and $M=800$, the program
gives

\begin{center}
\begin{tabular}{r|ccc}
        & $p_0=1/3$ & $p_0=1/2$ & $p_0=2/3$ \\
\hline
$t_0=6$ & 0.2032  & 0.2249 & 0.2489 \\
$t_0=12$& 0.2629  & 0.3973 & 0.4066 \\
$t_0=18$& 0.4019  & 0.5222 & 0.5535 \\
$t_0=24$& 0.4485  & 0.6173 & 0.6630 \\
$t_0=30$& 0.4838  & 0.6818 & 0.7454 \\
$t_0=36$& 0.4637  & 0.7306 & 0.8063 \\
$t_0=42$& 0.5271  & 0.7682 & 0.8578 \\
$t_0=48$& 0.5448  & 0.7978 & 0.8859 \\
\hline
\end{tabular}
\end{center}

\vspace{3mm}
Such values allow to find initial values $t_0,\alpha_0$ in order to have an wide probability to get
infinite convergent sequences $´\{p_n\}$ and $\{X_n\}$.

\section{ Appendix: Matlab programs}

{\bf Function   random\_value } \vspace{-4mm}
\begin{verbatim}

function [ak] = random_value(a,p)     %%% choose a random value in a with
                                      %%% probabilities proportional to p
  x=rand;
  ak=a(1);
  pp=cumsum(p/sum(p));
  for k=1:size(p,2)
      if x>pp(k)
          ak=a(k+1);
      end
  end
\end{verbatim}

{\bf Function   p\_sequence } \vspace{-4mm}
\begin{verbatim}

function [sec_p,tau]=p_sequence(m,alfa,beta,A,B,u,r,v,s)
                    %%% simulates m steps in an urn
                    %%% scheme with given parameters
    unif=rand(1,m); % random numbers in (0,1)
    t=zeros(1,m);   % total number of balls
    az=zeros(1,m);  % number of amber balls
    p=zeros(1,m);   % poportion of amber balls
    y=zeros(1,m);   % color of extracted ball
    az(1)=alfa;
    t(1)=alfa+beta;
    p(1)=az(1)/t(1);
    tau=0;
    k=1;
    while k<m && p(k)<=1 && p(k)>=0
       y(k)=(unif(k)<p(k));  % 1 if ineq. holds, 0 otherwise
       a=random_value(u,r);
       b=random_value(v,s);
       az(k+1)=az(k)+y(k)*(A-a)+(1-y(k))*b;
       t(k+1)=t(k)+y(k)*A+(1-y(k))*B;
       p(k+1)=az(k+1)/t(k+1);
       k=k+1;
    end
    if k<m && p(k)>1
       tau=k; p(k+1:m)=ones(1,m-k);
    elseif k<m && p(k)<0
        tau=k;
    else
        tau=0;
    end
    sec_p=p;
end

\end{verbatim}

{\bf Simulation of n p\-sequences of length m}\vspace{-2mm}

\begin{verbatim}

n=10;       % number of trajectories
m=5000;     % length of trajectories
alfa=30;    % inicial number of amber balls
beta=30;    % inicial number of blue balls
A=7; B=2; D=A-B;  % total number of added balls
u=[-5,-2,4,7]; r=[1,2,2,1];  % distribution of added blue balls
                             % when an amber ball is extracted
v=[-5,0,4]; s=[2,3,1];       % distribution of added amber balls
                             % when a blue ball is extracted
vtau=zeros(1,n);  % end of each trajectory
for iter=1:n
    [p,tau]=p_sequence(m,alfa,beta,A,B,u,r,v,s);
    plot(1:m,p)   % plot of each trajectory
    hold on
    vtau(iter)=tau;
end
tabulate(vtau)  % distribution of vtau
am=u*r'/sum(r)  % mean number of added blue balls
bm=v*s'/sum(s)  %  mean number of added amber balls

%%% plot of limit lines
if A==B & (am+bm>0 | am+bm<0)
   pstar=bm/(am+bm)
   plot([1,m],[pstar,pstar],'r','LineWidth',2)
elseif am==0 &bm>=0
      pstar=1
      plot([1,m],[pstar,pstar],'r','LineWidth',2)
elseif am>0 & bm>0
      pstar= (D-am-bm+((D-am-bm)^2+4*D*bm)^(1/2))/(2*D)
      plot([1,m],[pstar,pstar],'r','LineWidth',2)
   elseif am<0 & bm<0
      pstar= (D-am-bm-((D-am-bm)^2+4*D*bm)^(1/2))/(2*D)
      plot([1,m],[pstar,pstar],'r','LineWidth',2)
   elseif am>0 & bm==0 & ( (abs(am+bm)<D & (D-am-bm)^2+4*D*bm<0) | abs(am+bm)>D)
      pstar=0
      plot([1,m],[pstar,pstar],'r','LineWidth',2)
   elseif am>0 & bm<=0 & abs(am+bm)<D & (D-am-bm)^2+4*D*bm>0
     pstar1= (D-am-bm-((D-am-bm)^2+4*D*bm)^(1/2))/(2*D)
     pstar2= (D-am-bm+((D-am-bm)^2+4*D*bm)^(1/2))/(2*D)
     plot([1,m],[pstar1,pstar1],'r','LineWidth',1)
     plot([1,m],[pstar2,pstar2],'r','LineWidth',2)
end
\end{verbatim}

{\bf Estimation of $\P\{\tau=\infty\}$}\vspace{-4mm}
\begin{verbatim}

t1=6;    % minimal initial number of balls
t2=48;   % maximal initial number of balls
A=7; B=2;
u=[-5,-2,4,7]; r=[1,2,2,1];
v=[-5,0,4]; s=[2,3,1];
M=800;
q1=ones(t2+M*A,t2+M*A+1);
n=M-1
while n>=0
   q2=zeros(t2+n*A,t2+n*A+1);
   for t=t1+n*B:t2+n*A
       for a=1:t+1
           x=0; y=0;
           for k=1:size(u,2)
              if a-1+A-u(k)>=0 & a-1-u(k)<=t
                 x=x+q1(t+A,a+A-u(k))*r(k)/sum(r);
              end
           end
           for k=1:size(v,2)
              if  a-1+v(k)>=0 & a-1+v(k)<=t+B
                 y=y+q1(t+B,a+v(k))*s(k)/sum(s);
              end
            end
           q2(t,a)=x*(a-1)/t+y*(1-(a-1)/t);
       end
   end
   %pause
   q1=q2;
   n=n-1
end
\end{verbatim}

\end{document}